\documentclass[12pt]{amsart}
\usepackage[dvips]{graphicx}
\usepackage{amsmath, amsthm, amsfonts, latexsym, amscd, amssymb, 
}
\newtheorem{thm}{Theorem}[section]
\newtheorem{cor}[thm]{Corollary}
\newtheorem{df}[thm]{Definition}
\newtheorem{lem}[thm]{Lemma}
\newtheorem{prop}[thm]{Proposition}

\newtheorem{ex}[thm]{Example}
\numberwithin{equation}{section}
\def\proof{\textsc{Proof:}\ }
\def\endproof{$\Box$ \medskip}

\begin{document}
\title{Public key exchange using right transversals and  right loops}
\author{Akhilesh Chandra Yadav$^1$ and Vipul Kakkar$^2$}
\address{$^1$Department of Mathematics\\
M G Kashi Vidyapith, Varanasi, India}
\email[A.~C.~Yadav]{akhileshyadav538@gmail.com }
\address{$^2$School of Mathematics\\
Harish-Chandra Research Institute, Allahabad, India}
\email[V.~Kakkar]{vplkakkar@gmail.com}

\date{}
\maketitle
\begin{abstract}
In this article, we describe  a key exchange protocol  based on right transversals. We also describe it for general extension associated to right loops.
\end{abstract}


\section{Introduction}
Loops as algebraic structures had been an important object of study in mathematics starting from the first half of the 20th century with the works of Baer~\cite{baer}, Albert~\cite{albert1, albert2} and Bruck~\cite{bruck}. 
It has always been practice to study right loops through its  right inner mappings and  inner mapping groups (also called group torsions~\cite{lal}). 
The notion of general extension and the notion of general extension associated to a right loop has been given\cite{lal}. Indeed, it is observed that every right loop $S$ can be embedded into a group as  right transversal with some universal property\cite{lal}. The smallest subgroup generated by $S$ is the group $G_S\times S$, the general extension associated to $S$. 
The notion of encryption and decryption has been introduced for  a given right loop\cite{yadav} by using bracket arrangement of weight $n$.

Diffie-Hellman scheme\cite{lpilz} is a key exchange system for establishing a common key between $A$ and $B$. Their key exchange protocol  is based on Cyclic groups. In 2010, M. Habeeb, D. Kahrobae, C. Koupparie and V. Spilrain introduced a new concept  of key exchange protocol using semi-direct product  of (semi) groups and then focused on practical instances of this general idea \cite{hkkv}. Their concept motivates us to give key exchange protocol based on a right transversal to a subgroup of the given group (general extensions) and general extensions  associated to right loops. We also discuss key exchange protocol in the case of right gyrogroups\cite{ylal} and twisted right gyrogroups\cite{ylal1}.

\section{c-groupoids and general extensions}

A group $G$ is called  a general extension of a group $H$ if $H$ can be treated as a subgroup of $G$. Let $H$ be a subgroup of group $G$.  A {\it right transversal} to $H$ in $G$ is a subset $S$ of $G$ obtained by
selecting one and only one member from each right coset of $G$ mod $H$ including identity $e$ of group $G$. It is observed that each right transversal $S$ to a subgroup $H$ in a group $G$ determines an algebraic structure $(S, H, \sigma, f)$\cite{lal} in the sense of following:  
\begin{df} 
A quadruple $(S,H,\sigma ,f)$, where $S$ is a groupoid with identity  $e$, $H$ a group which acts on $S$ from right through  a given action  
$\theta$, $\sigma$ a map from
$S$ to $H^H$ (the set of all maps
 from  $H$ to $H$) and $f$ a map from 
$S\times S$ to $H$, is called a  c-groupoid if it satisfies the following 
conditions:
\begin{enumerate}
\item $ x\circ y\ =\ y\; \Rightarrow x\ =\ e$,
\item   For each $x\in S,\exists x^{\prime}\in S$ such that  
$x^{\prime}\circ x\ =\ e$,
\item $\sigma_{e}\ =\ I_{H}$, the identity map on $H$, where  $\sigma_{x}$ denotes the image of $x$
  under the map $\sigma$,
\item $ f(x,e)\ =\ f(e,x)\ =\ 1$, the identity  of $H$,
\item $\sigma_{x}(h_{1}h_{2})\ =
\ \sigma_{x}(h_{1})\sigma_{x\theta h_{1}}(h_{2})$,
\item $(x\circ y)\circ z\ =\ x\theta f(y,z)\circ (y\circ z)$,
\item $(x\circ y)\theta h\ =\ x\theta \sigma_{y}(h)\circ 
(y\theta h)$,
\item $f(x,y)f(x\circ y,z)\ =\ \sigma_{x}(f(y,z))f(x\theta f(y,z), y\circ z)$,
\item $f(x,y)\sigma_{x\circ y}(h)\ =
\ \sigma_{x}(\sigma_{y}(h))f(x\theta
\sigma_{y}(h),y\theta h)$.
\end{enumerate}
where $x,y,z\in S$ and $h_{1},h_{2},h\in H$. 
\end{df}
Conversely, we have
\begin{thm}(\cite{lal}, Theorem 2.2)  \label{cgth1} 
Given a  $c$-groupoid $(S, H,\sigma, f)$ there is a group $G$ which contains 
$H$ as a  subgroup and $S$ as a right transversal of $H$ in $G$ such that the corresponding $c$-groupoid is  $(S,H,\sigma, f)$.
\end{thm}

The corresponding group $G$ is $H\times S$ together with a binary operation $.$ given by :
\begin{eqnarray}\label{gext1}
 (a, x) \cdot (b, y) & = & \left( a \sigma_{x}(b) f(x \theta b, y),  (x \theta b) \circ y \right). 
\end{eqnarray}
This group is termed as general extension of a group $H$ by a set $S$ satisfying above properties (or a general extension associated to a c-groupoid).

Let $(S, H, \sigma, f)$ be a c-groupoid. Let $x\in S$ and $a\in H$. 
Define $\beta^n (x,a)$ inductively by the following:
$$
\beta^1 (x,a)= x,\quad \beta^2 (x,a)= x\theta a\; o\; x. 
$$
If $\beta^r (x,a)$ is defined then 
\begin{eqnarray}\label{bn}
\beta^{r+1} (x,a)= \beta^r (x,a)\theta a \; o\;  x.
\end{eqnarray}
 Similarly, we define $g^n (x, a)$inductively by the following:
$$g^1(x, a)= a, \qquad g^2(x, a)= a\sigma_x(a) f(x\theta a, x).$$
 If $g^n (x, a)$ is defined, then
 \begin{eqnarray}\label{gn}
g^{n+1} (x, a)= g^n(x,a) \sigma_{\beta^n(x, a)} (a) f(\beta^n(x, a)\theta a, x)
\end{eqnarray}

Using Eqs. \ref{gext1},\ref{bn}, \ref{gn}, we have
\begin{equation}\label{cyclic}
(a, x)^n = (g^n(x, a), \beta^n (x, a))
\end{equation} 
for each $n\in \mathbb{N}$. 

Since 
\begin{eqnarray*}
(a, x)^n .(a, x)^m= (a, x)^{m+n}= (a, x)^m . (a, x)^n
\end{eqnarray*}
Using Eqns. \ref{cyclic} and \ref{gext1}, we have the following:

\begin{lem}\label{loop}
Let $(S, H, \sigma, f)$ be a c-groupoid and $x\in S\setminus\{e\}, h\in H\setminus \{1\}$. Then
\begin{enumerate}
\item \begin{eqnarray*}
 g^{n+m}(x, a)& = & g^n(x, a)\sigma_{\beta^n(x, a)} (g^m(x, a))  f(\beta^n(x, a)\theta g^m(x, a), \beta^m(x, a)) \\
& = & g^m(x, a)\sigma_{\beta^m(x, a)} (g^n(x, a)) f(\beta^m(x, a)\theta g^n(x, a), \beta^n(x, a)) 
\end{eqnarray*}

\item $\beta^m(x, a)\theta g^n(x, a)o \beta^n(x, a)= \beta^n(x, a)\theta g^m(x, a)o \beta^m(x, a)=  \beta^{m+ n}(x, a)$.
\end{enumerate}
\end{lem}
Next, for $x\in S, a\in H$, we define $[a \sigma_x (a)]_m$ inductively by :
$$
[a \sigma_x (a)]_0= a,\qquad [a \sigma_x (a)]_1= a \sigma_x (a)
\quad {\rm and}\;  
[a \sigma_x (a)]_n= a\sigma_x ([a \sigma_x (a)]_{n-1}).
$$

Then, we have the following:
\begin{lem}\label{bracket} Let $(S, H, \sigma, f)$ be a c-groupoid and $x\in S\setminus\{e\}, h\in H\setminus \{1\}$. For $m\ge 2$,  $\beta^m (x, a)$ is given by 
\begin{eqnarray}
 \left(\left(\ldots (x\theta [a\sigma_x (a)]_{m-2}\ o \ x\theta [a\sigma_x (a)]_{m-3})\ o \ldots\right)\ o\ x\theta [a\sigma_x (a)]_{0} \right)o x &  &
\end{eqnarray}\end{lem}

\proof
For $m\ge 2$,
\begin{eqnarray*}
\beta^m (x, a) &= & \{\beta^{m-1}(x, a)\theta a\} o x\\
& = &\{ (\{\beta^{m-2}(x, a)\theta a\} o x)\theta a\} o x\\
& = &  ({\beta^{m-2}(x, a)\theta (a\sigma_x (a))} o {x\theta a}) o x\\
&  = & (\beta^{m-2}(x, a)\theta {[a\sigma_x (a)]_1 }\; o \; x\theta {[a\sigma_x (a)]_0}) \; o\;  x\\
& = & (((\beta^{m-3}(x, a)\theta a\;  o\;  x)\theta {[a\sigma_x (a)]_1})\; o\; {x\theta [a\sigma_x (a)]_0}) \; o\; x\\
& = &  ((\beta^{m-3}(x, a)\theta (a\sigma_x ([a\sigma_x (a)]_1))\;  o\;  x\theta [a\sigma_x (a)]_1)\;  o\;  x\theta [a\sigma_x (a)]_0 ) o x\\
& = &  ((\beta^{m-3}(x, a)\theta ([a\sigma_x (a)]_2)\;  o\;  x\theta [a\sigma_x (a)]_1)\;  o\;  x\theta [a\sigma_x (a)]_0 ) o x\\
& = & \ldots\qquad \ldots \qquad \ldots \\
& = & \ldots\qquad \ldots \qquad \ldots \\
& = &  [[\ldots [\beta^1 (x, a)\theta ([a\sigma_x (a)]_{m-2})\;  o\;  x\theta [a\sigma_x (a)]_{m-3}]\;  o \;\ldots] \; o\;  x\theta [a\sigma_x (a)]_0 ]\;  o\; x\\
& = &  [[\ldots [x\theta ([a\sigma_x (a)]_{m-2})\;  o\;  x\theta [a\sigma_x (a)]_{m-3}]\;  o \;\ldots] \; o\;  x\theta [a\sigma_x (a)]_0 ]\;  o\; x
\end{eqnarray*}
\endproof
 \begin{cor}\label{sigma}
 If $\sigma_x = I_H$, then $[a\sigma_x(a)]_m= a^{m+1}$ and so 
\begin{eqnarray}
\beta^m (x, a)& = & \left(\left(\ldots (x\theta a^{m-1}\ o \ x\theta a^{m-2})\ o \ldots\right)\ o\ x\theta a \right)o x.
\end{eqnarray}
\end{cor}

\begin{cor}\label{eta}
 If $\sigma_x (h) = \eta(h)$ for all $x\in S\setminus\{e\}$ and $h\in H$, where $\eta\in Aut\ H$ is an involution. Then 
\[ [a \sigma_x (a)]_m = \left\{\begin{array}{ccc}
a[\eta(a)a]^{\frac{m}{2}} & {\rm if }\  m \ {\rm is\, even}\\

& & \\

[a\eta(a)]^{\frac{m+1}{2}} & \rm{ if}\  m\  \rm{is\, odd}
\end{array}\right. \]
Thus 
\begin{eqnarray*}\label{taut}
\beta^m (x, a)& = & \left\{\begin{array}{lcr}
\left(\left(\ldots (x\theta [a (\eta(a) a)^\frac{m-2}{2}]\ o \ x\theta [a\eta(a)]^\frac{m-2}{2})\ o \ldots\right)\ o\ x\theta a \right)\, o x & \rm{if}\  m\ \rm{is\ even}\\

 &   & \\
 \left(\left(\ldots (x\theta [a\eta(a)]^\frac{m-1}{2}\ o \ x\theta [a(a\eta(a)^\frac{m-3}{2})\ o \ldots\right)\ o\ x\theta  a \right)\ o x & \rm{if}\, m\, \rm{is\ odd}.
\end{array}\right.
\end{eqnarray*}
\end{cor}

\section{Key exchange protocol using right transversals}
Let $S$ be a right transversal to a subgroup $H$ in a group $G$. Clearly each $g\in G$ can be uniquely expressed as $g= hx$ for $h\in H$ and $x\in S$. For the sake of convenience, we shall call $h$ as {\it subgroup component}  of $g$ and $x$ as {\it representative} of $g$. Let $x\in S\setminus \{e\}$ and $a\in H\setminus \{1\}$. Keeping $x$ and $a$ as public, Alice chooses her private key $m\in \mathbb{N}$ and Bob chooses his private key $n\in \mathbb{N}$. Both are agree to work with the cyclic subgroup $\{(a x)^r= g^r(x,a) \beta^r (x, a); |\; r\in\mathbf{N}\}\cup\{1 e\}$. In this case the  key exchange  protocol is given by:
\begin{enumerate}
\item Alice computes $(a x)^m=g^m(x,a)\beta^m (x, a)$ and sends only its  representative  $\beta^m (x, a)$ (given by \ref{bracket})  to the Bob.
\item Bob computes $(a x)^n= g^n(x,a) \beta^n (x, a)$ and sends only its representative  $\beta^n(x, a)$ (given by \ref{bracket})  to the Alice.
\item Alice computes 
$
(\phi \beta^n (x, a) ).(g^m(x,a) \beta^m (x, a))$ which is equal to 
$$
(\phi \sigma_{\beta^n (x, a)} f(\beta^n (x, a) \theta g^m(x,a), \beta^m (x, a)))\; [\beta^n (x, a) \theta g^m(x,a)\; o\; \beta^m (x, a) ] 
$$
and her key $K_A=  \beta^n (x, a) \theta g^m(x,a)\; o\; \beta^m (x, a)= \beta^{m+n} (x, a)$. Note that Alice is not able to compute subgroup component because she does not know the subgroup component $\phi= g^n(x, a)$. Thus, she computes only the representative.
\item Bob computes 
$
(\psi \beta^m (x, a) ).(g^n(x,a) \beta^n (x, a))$ which is equal to
$$
 (\psi \sigma_{\beta^m (x, a)}(g^n(x,a)) f(\beta^m (x, a) \theta g^n(x,a), \beta^n (x, a)))  [\beta^m (x, a) \theta g^n(x,a)\; o\; \beta^m (x, a)] 
$$
and his key $K_B=  \beta^m (x, a) \theta g^n(x,a)\; o\; \beta^n (x, a)= \beta^{n+m} (x, a)$. Note that Bob is not able to compute subgroup component because he does not know $\psi =g^m(x, a)$. Thus, he computes only the representative.
\item
In the general extension $G= H\times S$ determined by a c-groupoid $(S, H, \sigma, f)$,
$$
(\phi, \beta^n).(\psi, \beta^m)= (\psi, \beta^m). (\phi, \beta^n)= (a,x)^{n+m}.
$$
Thus the shared common key is $K= K_A= K_B$.
\end{enumerate}
\section{Right loops and Key exchange protocol}

A non empty set $S$ together with binary operations $o$ is called a {\it right loop} if for each $x, y\in S$, the equation $Xo x = y$ has a unique solution in $S$. The identity element  of $S$ is denoted by $e$.

Let $(S, o)$ be a right loop with identity $e$ and $y, z$ in $S$. 
 The map $f(y, z)$ from $S$ to $S$ given by the equation 
\begin{eqnarray}\label{bas1}
f(y,z)(x)o (yoz) & = & (xoy)oz,\ \qquad x\in S
\end{eqnarray}
belongs to $Sym\ S$ (the Symmetric group on $S$) and is called a right inner mapping of $(S, o)$.
Indeed $f(y, z)\in Sym (S\setminus\{e\})\subseteq Sym\ S$.
The subgroup $G_S$ of $Sym (S\setminus\{e\})\subseteq Sym\ S$ generated by $\{f(y, z)\ | \, y, z\in S\}$ is  called the right inner mapping group (also called  the group torsion \cite{lal}) of $(S, o)$.

Further, let $h\in Sym (S\setminus\{e\})\subseteq Sym\ S$ and $y\in S$.  
Define $\sigma_y (h)\in Sym (S\setminus\{e\})\subseteq Sym\ S$ by the equation
\begin{eqnarray}\label{rlal}
h(x o y)  & = &  \sigma_y (h) (x) o h(y),  \  \qquad x\in S
\end{eqnarray}

 For the sake of convenience  we shall also write $x\theta h$ for $h(x)$. Thus the equations (\ref{bas1}) and (\ref{rlal}) also read as 
\begin{eqnarray}\label{bas2}
x\theta f(y,z)(x)o (yoz) & = & (xoy)oz,\ \qquad x\in S
\end{eqnarray}
and 
\begin{eqnarray}\label{rlal2}
(xoy)\theta h & =  & x\theta \sigma_y (h) o y\theta h
\end{eqnarray} 
respectively.

 \begin{prop}\label{rlal1}
Let $(S, o)$ be a right loop with identity $e$. Then it determines a c-groupoid $(S, G_S, \sigma, f)$\cite{lal}. 
\end{prop}

The group $G_S\times S$ determined by c-groupoid $(S, G_S, \sigma, f)$ is the smallest group generated by $S$, which is known as the general extension associated to the right loop $S$. Since all the results described as above hold in the general extension $G_S\times S$, therefore Bob and Alice may agree to work with the given right loop. Their key exchange protocol is described by the following:

\indent Let $S$ be a right loop and $x\in S\setminus \{e\}$ and $a\in G_S\setminus \{I_S\}$. Keeping $x$ and $a$ as public, Alice chooses her private key $m\in \mathbb{N}$ and Bob chooses his private key $n\in \mathbb{N}$. Both are agree to work with the cyclic subgroup $\{(a, x)^r= (g^r(x,a), \beta^r (x, a)); |\; r\in\mathbf{N}\}\cup\{(I_S,  e)\}$. In this case the  key exchange  protocol is given by:
\begin{enumerate}
\item Alice computes $(a, x)^m= (g^m(x, a), \beta^m (x, a))$ and sends only the second component $\beta^m (x, a)$ to the Bob.
\item Bob computes $(a, x)^n= (g^n(x, a), \beta^n (x, a))$ and sends only the second component $\beta^n (x, a)$ to  the Alice.
\item Alice computes 
$ \beta^n (x, a) \theta g^m(x,a)\; o\; \beta^m (x, a) ) $. Her key is now $$K_A=  \beta^n (x, a) \theta g^m(x,a)\; o\; \beta^m (x, a).$$
\item Bob computes 
$ \beta^m (x, a) \theta g^n(x,a)\; o\; \beta^m (x, a) $. His key is now  $$K_B=  \beta^m (x, a) \theta g^n(x,a)\; o\; \beta^n (x, a).$$
\item  Using Lemma \ref{loop}, 
 the shared  common key $K= K_A= K_B$.
\end{enumerate}

Using Lemma (\ref{sigma}), we have: 
\begin{cor}
If $(S, o)$ is a right gyrogroup\cite{ylal}, then  the shared common key will be $\beta^{m+n} (x, a) $, where  $\beta^m (x, a)= [(\ldots((x\theta a^{m-1}\ o\ x\theta a^{m-2}) o x\theta a^{m-3})\ldots)o x\theta a] o x $ for $m\in \mathbb{N}$.  
\end{cor}

\begin{cor}
If $(S, o)$ is a twisted right gyrogroup\cite{ylal1}, then  $\sigma_x $ for every $x\in S\setminus \{e\}$, is a fixed involutory automorphism of $G_S$  say $\eta$ . In this case the shared common key will be $\beta^{n+m} (x, a)$, where $\beta^n (x, a)$ is given by Corollary \ref{eta}.
\end{cor}

\begin{ex}
Let $S = \{e,\ x_1,\  x_2,\ \ldots, \ x_{15}\}$. Define a binary operation $o$ on $S$ by taking $e$ as the identity and defining $x_i \ o \ x_j = x_i$ if $i\ne j$ and $x_i\ o \ x_i = e$. Then $(S, o)$ is a right loop with $x' = x$ and $f(x', x) = I_S$ for all $x\in S$. Also, for $i\ne j$
\begin{displaymath}
x_k\theta f(x_i, x_j)  =  \left\{ 
\begin{array}{lcr}
x_k &{\rm for} & i, j\neq k\\
x_j & {\rm for} & k=i\\
x_i & {\rm for} & k=j
\end{array}
\right.
\end{displaymath}
This shows that the group torsion $G_S = Sym\ (S\setminus \{e\})$. It is also evident that $f(x_i, x_j)\in Aut\ (S, o)$. Thus, $Aut\ (S, o) = G_S$ and so $(S, o)$ is a right gyrogroup~\cite{ylal}. 
Take $x= x_3$ and $a= (x_3\ x_4\  x_1\  x_9\  x_8\  x_7)$. Then 

\begin{eqnarray*}
 \beta^1 (x, a)= x  & & g^1(x, a)= a\\
 \beta^2 (x, a)= x_4 o x_3= x_4& & g^2(x, a)= a^2 (x_3\  x_4)=(x_3\  x_1\  x_8\  x_4\  x_9\  x_7)\\
\beta^3 (x, a) = x_9 o x_3= x_9 & &\\
 g^3(x, a) =&  (x_3\ x_1\  x_8\  x_4\  x_9\  x_7)\ a\ (x_1\  x_3) & = (x_1\  x_7\, x_4\, x_8\, x_3\, x_9)
\end{eqnarray*}
Alice chooses her private number $2$ and sends $\beta^2 (x, a)= x_4$ to Bob. Bob chooses his private number $3$ and sends $\beta^3(x, a)= x_1$. Now, Alice computes 
$$x_1\theta g^2 (x, a) o x_4=x_1\theta (x_3\  x_1\  x_8\  x_4\  x_9\  x_7) o x_4= x_8 $$
 and Bob computes 
$$x_4\theta g^3 (x, a) o x_1=x_4\theta (x_1\  x_7\, x_4\, x_8\, x_3\, x_9) o x_1= x_8 o x_1= x_8 .$$
 Thus their shared common key is $x_8$. 
\end{ex}


\begin{thebibliography}{breitestes Label}
\bibitem[1]{albert1} A. A. Albert, {\it Quasigroups I}, Trans. Amer. Math. Soc. {\bf 54} (1943), 507-519. 
\bibitem[2]{albert2} A. A. Albert, {\it Quasigroups II}, Trans. Amer. Math. Soc. {\bf 55} (1944), 401 -419.

\bibitem[3]{baer} R. Baer, {\it Nets and groups}, Trans. Amer. Math. Soc. {\bf 46} (1939), 110- 141.
\bibitem[4]{bruck} R. H. Bruck, {\it Contributions to the theory of Loops}, Trans. Amer. Math. Soc. {\bf 60} (1946), 245- 354.
\bibitem[5]{hkkv} M. Habeeb,  D.  Kahrobaei, C. Koupparis and V. Shpilrain, {\it  Public key exchange using semi-direct product of (semi)groups}, preprint available at http://arxiv.org/abs/1304.6572. 

\bibitem[6]{lal} R. Lal, {\it Transversals in groups}, Journal of algebra, {\bf 181}(1996)\  70-81.
\bibitem[7]{lpilz} R. Lidl and G. Pilz, {\it Applied Abstract Algebra}, Second edition, Springer (First Indian Reprint, 2004)
\bibitem[8]{ylal} R. Lal and A. C. Yadav, {\it Topological right gyrogroups and  gyrotransversals},  Communications in Algebra, {\bf  41(09)}(2013)\  3559 - 3575.
\bibitem[9]{ylal1}  R. Lal and A. C. Yadav, {\it Twisted Automorphisms and Twisted right Gyrogroups}, accepted in Communications in Algebra, estimated publication date 18 Dec, 2014(online).
\bibitem[10]{yadav} A. C. Yadav,{\it Generating non-isomorphic right loops of a given order}, J. Disc. Math. Sci. and Cryptography, Vol. 16 (2, 3), 139-148, 2013.
\end{thebibliography}
\end{document}